\documentclass{article}
\usepackage{amsmath}
\usepackage{amssymb}
\newtheorem{theorem}{Theorem}

\title{Comparison of two techniques for\\
proving nonexistence of strongly regular graphs}
\author{Va\v{s}ek Chv\'{a}tal\thanks{Department of Computer Science and Software Engineering, Concordia University, Montr\'{e}al, Qu\'{e}bec, Canada}}
\date{}

\begin{document}

\maketitle
\begin{abstract}
  We show that the method of counting closed walks in strongly regular
  graphs rules out no parameter sets other than those ruled out by the
  method of counting eigenvalue multiplicities.
\end{abstract}

Following Bose \cite{Bo}, a {\em strongly regular graph with
  parameters $n,k,\lambda,\mu$\/} means an undirected graph $G$ such
that
\begin{itemize}
\item $G$ has $n$ vertices,
\item $G$ is regular of degree $k$,
\item every two adjacent vertices of $G$ have precisely $\lambda$ common
neighbours, 
\item every two nonadjacent vertices of $G$ have precisely $\mu$
common neighbours.
\end{itemize}
Complete graphs have these four properties (with $k=n-1$,
$\lambda=n-2$, and any $\mu$) and so have their complements (with
$k=n-1$, any $\lambda$, and $\mu=n-2$). Let us follow the convention
of excluding these trivial examples from the class of strongly regular
graphs: let us assume that 
\begin{equation}\label{bounds}
0\;<\;k\;<\;n-1.
\end{equation}
If there exists a strongly regular graph with parameters
$n,k,\lambda,\mu$, then 
\begin{equation}\label{srg1}
(n-1-k)\mu \;=\; k(k-1-\lambda ).
\end{equation}
(This identity follows directly from counting in two different ways
all sequences $w_0,w_1,w_2$ of vertices $w_0,w_1,w_2$ such that $w_0$
is prescribed, $w_0,w_1$ are nonadjacent, $w_1,w_2$ are adjacent, and
$w_0,w_2$ are nonadjacent: choosing $w_2$ first and $w_1$ second gives
the left-hand side; choosing $w_1$ first and $w_2$ second gives the
right-hand side). Another widely known condition that is necessary for
the existence of a strongly regular graph with parameters
$n,k,\lambda,\mu$ goes as follows:
\begin{theorem}\label{eigen}
  If there exists a strongly regular graph with parameters
  $n,k,\lambda,\mu$, then 
\begin{equation}\label{srg2}
\frac{1}{2}\left(n-1 \pm \frac{(n-1)(\lambda-\mu)+2k}{\sqrt{(\lambda-\mu)^2+4(k-\mu)}}\right)
\;\;\mbox{ are nonnegative integers. }
\end{equation}

\end{theorem}
The method used in the proof of Theorem \ref{eigen} can be traced back
to Connor and Clatworthy \cite{CC}; it was used by Hoffman and
Singleton \cite{HS} in the special case $\lambda=0$, $\mu=1$ and by
Wilf \cite{W} in the special case $\lambda=\mu=1$. Additional
information on strongly regular graphs can be found in \cite{BvL,C}
and elsewhere.\\

The famous Friendship Theorem of Erd\H{o}s, R\'{e}nyi, and S\'{o}s
(\cite{ERS}, Theorem 6) states that
\begin{quote}
{\em if, in a finite undirected graph $G$, every two vertices have
precisely one common neighbour, then some vertex of $G$ is adjacent to
all the vertices of $G$ except itself.}
\end{quote}
It is relatively easy to show that every counterexample $G$ to this
theorem would have to be regular (the first step is proving that every
two nonadjacent vertices must have the same degree).  The rest of the
proof amounts to proving that there is no strongly regular graph with
parameters $n,k,1,1$. For this purpose, Erd\H{o}s, R\'{e}nyi, and
S\'{o}s invoke a theorem of Baer \cite{B}, whose special case asserts
that every polarity in a projective plane of order at least $2$ maps
some point to a line that contains this point. To make this proof of
the Friendship Theorem self-contained, one may extract from \cite{B}
the corresponding fragment of Baer's reasoning; this is precisely what
Longyear and Parsons~\cite{LP}, and later also Huneke~\cite{H}, seem
to have done.  The argument generalizes to the context of strongly
regular graphs as follows.
\begin{theorem}\label{walks}
  If there exists a strongly regular graph with parameters
  $n,k,\lambda,\mu$, then every prime $p$ divides the integer $c_p$ 
  defined by the recurrence 
\begin{equation}\label{rec}
c_{\ell} \,=\, \mu nk^{\ell-2} + (\lambda-\mu) c_{\ell-1} + (k-\mu)c_{\ell-2} 
\end{equation}
with the initial conditions $c_{0} = n$, $c_{1} = 0$.
\end{theorem}
\noindent{\bf Proof.}  A {\em walk of length $\ell$\/} in a graph $G$ is a
sequence $w_0,w_1,\ldots ,w_\ell$ of (not necessarily distinct)
vertices such that each $w_i$ with $0\le i <\ell$ is adjacent to
$w_{i+1}$.  The walk is called {\em closed\/} if $w_\ell=w_0$ and it
is called {\em open\/} if $w_\ell\ne w_0$. If there is a strongly
regular graph $G$ with parameters $n,k,\lambda,\mu$, then the number
of closed walks of length $\ell$ in $G$ satisfies the recurrence for
$c_\ell$, since 

\begin{itemize}
\item $k c_{\ell-2}$ of these walks have $w_{\ell-2}=w_{0}$,
\item $\lambda c_{\ell-1}$ of these walks have $w_{\ell-2}\ne w_{0}$ with $w_{\ell-2}$, $w_{0}$ adjacent, 
\item $\mu(nk^{\ell-2} - c_{\ell-2} - c_{\ell-1})$ of these walks have $w_{\ell-2}\ne w_{0}$ 
with $w_{\ell-2}$, $w_{0}$ nonadjacent.
\end{itemize}
The proof is completed by observing that every prime $p$ divides the
number of closed walks of length $p$ in $G$, since each equivalence class of
the equivalence relation $\sim$ defined on the set of all closed walks
of length $p$ in $G$ by
setting 
\begin{quote}

  $u_0,u_1,\ldots ,u_{p-1},u_0 \;\;\sim\;\; v_0,v_1,\ldots ,v_{p-1},v_0$ if and only if\\
  there is an integer $s$ such that $0\le s <p$ and $v_i=u_{(i+s)
    \bmod p}$ for all $i=0,1,\ldots ,p-1$ 
\end{quote}
consists of $p$ distinct walks. \hfill $\Box$\\

The purpose of this note is to show that the necessary condition of
Theorem~\ref{walks} is subsumed in the necessary condition of
Theorem~\ref{eigen}:
\begin{theorem}\label{comp}
  Let $n,k,\lambda,\mu$ be nonnegative integers with properties
  {\rm(\ref{bounds})}, {\rm(\ref{srg1})}, {\rm(\ref{srg2})}. Then
  every prime $p$ divides the integer $c_p$ defined by the recurrence
  {\rm (\ref{rec})} with the initial conditions $c_{0} = n$, $c_{1} =
  0$.
\end{theorem}

\noindent{\bf Proof.} We have 
\[
c_\ell = k^\ell+m_1r_1^\ell+m_2r_2^\ell
\]
with
\begin{eqnarray*}
m_1 &=& \frac{1}{2}\left(n-1 - \frac{(n-1)(\lambda-\mu)+2k}{\sqrt{(\lambda-\mu)^2+4(k-\mu)}}\right),  \\ 
m_2 &=& \frac{1}{2}\left(n-1 + \frac{(n-1)(\lambda-\mu)+2k}{\sqrt{(\lambda-\mu)^2+4(k-\mu)}}\right),\\
r_1 &=& \frac{1}{2}\left((\lambda-\mu) + \sqrt{(\lambda-\mu)^2+4(k-\mu)}\right),  \\ 
r_2 &=& \frac{1}{2}\left((\lambda-\mu) - \sqrt{(\lambda-\mu)^2+4(k-\mu)}\right).  
\end{eqnarray*}

{\sc Case 1:} $(n-1)(\lambda-\mu)+2k\ne 0$.  In this case, property (\ref{srg2}) implies that 
$\sqrt{(\lambda-\mu)^2+4(k-\mu)}$ must be rational. Since the square root of an integer is rational 
only if it is an integer, it follows that $\sqrt{(\lambda-\mu)^2+4(k-\mu)}$ is an
integer.  Now
\[
\sqrt{(\lambda-\mu)^2+4(k-\mu)}\bmod 2 = 
((\lambda-\mu)^2+4(k-\mu))\bmod 2 = 
(\lambda-\mu) \bmod 2,
\]
and so $r_1$ and $r_2$ are integers.  By Fermat's Little Theorem, 
\[
(k^p+m_1r_1^p+m_2r_2^p) \bmod p \;=\; (k+m_1r_1+m_2r_2) \bmod p\,; 
\]
the right-hand side is zero since $k+m_1r_1+m_2r_2 =c_1=0$.\\

{\sc Case 2:} $(n-1)(\lambda-\mu)+2k= 0$.  In this case, property
(\ref{bounds}) implies that $\mu-\lambda =1$, and so $n=2k+1$; in
turn, property (\ref{srg1}) implies that $k=2\mu$. Now
\[
c_\ell \;=\; (2\mu )^\ell+2\mu (r_1^\ell+r_2^\ell) 
\]
with
\[
r_1 \,=\, \frac{-1 + \sqrt{4\mu +1}}{2},\;\;
r_2 \,=\, \frac{-1 - \sqrt{4\mu +1}}{2}.
\]
Expanding $r_1^\ell+r_2^\ell$ as 
\[
2^{-\ell}\sum_{j=0}^\ell\binom{\ell}{j}\left(\sqrt{4\mu +1}\right)^j(-1)^{\ell-j}\;+\;
2^{-\ell}\sum_{j=0}^\ell\binom{\ell}{j}\left(-\sqrt{4\mu +1}\right)^j(-1)^{\ell-j},
\]
we conclude that 
\[
c_\ell \;=\; (2\mu )^\ell+2\mu (r_1^\ell+r_2^\ell) \;=\; (2\mu )^\ell+4\mu \left(\frac{-1}{2}\right)^{\!\ell}\sum_{i=0}^{\lfloor \ell/2\rfloor}\binom{\ell}{2i}(4\mu +1)^{i}.
\]
In particular, $c_2=2\mu (4\mu +1)$. When $p$ is an odd prime, we have $c_p=(2\mu)^p-4\mu x_p$ with
\[
x_p= 2^{-p}\sum_{i=0}^{(p-1)/2}\binom{p}{2i}(4\mu +1)^{i};
\]
note that $x_p$ may not be an integer, but $2^px_p$ is one; since $p$
divides every $\binom{p}{j}$ with $0<j<p$, we have $2^p x_p \bmod p
=1$. By Fermat's Little Theorem, $2^{p-1} \bmod p=1$; since $4\mu x_p$
is an integer (it equals $(2\mu)^p-c_p$), it follows that
\[
4\mu x_p \bmod p = (2^{p-1}\cdot 4\mu x_p) \bmod p = (2\mu\cdot 2^{p}x_p) \bmod p
= 2\mu \bmod p.
\]
By Fermat's Little Theorem again, $(2\mu)^{p} \bmod p=2\mu \bmod p$, and so
$c_p\bmod p = 0$.
\hfill $\Box$


\begin{thebibliography}{99}

\bibitem{B}  R.~Baer, Polarities in finite projective
planes, {\em Bull. Amer. Math. Soc.\/} {\bf 52} (1946), 77--93.

\bibitem{Bo} R.C.~Bose, Strongly regular graphs, partial geometries
  and partially balanced designs, {\em Pacific J. Math.\/} {\bf 13} (1963), 389--419.

\bibitem{BvL} A.E.~Brouwer and J.H.~van Lint, Strongly regular graphs
  and partial geometries, in: {\em Enumeration and design
    (D.M.~Jackson and S.A.~Vanstone, eds.)\/}, Academic Press, Toronto,
  1984, pp. 85--122.

\bibitem{C} P.J.~Cameron, Strongly regular graphs, in: {\em Topics in
    algebraic graph theory (L.W.~Beineke and R.J.~Wilson, eds.\/},
    Cambridge University Press, Cambridge, 2004, pp. 203--221.

\bibitem{CC} W.S.~Connor and W.H.~Clatworthy, Some theorems for partially balanced designs, 
{\em Ann. Math. Stat.\/} {\bf 25} (1954), 100--112.

\bibitem{ERS} P.~Erd\H{o}s, A.~R\'{e}nyi, and V.T.~S\'{o}s, On a
problem of graph theory, {\em Studia Sci.Math. Hungar.\/} {\bf 1}
(1966), 51--57.

\bibitem{H} C.~Huneke, The Friendship Theorem, {\em
    Amer. Math. Monthly\/} {\bf 109} (2002), 192--194.

\bibitem{HS} A.J.~Hoffman and R.R.~Singleton, On Moore graphs with
  diameters $2$ and $3$, {\em IBM J. Res. Develop.\/} {\bf 4} (1960),
  497--504.

\bibitem{LP} J.Q.~Longyear and T.D.~Parsons, The friendship theorem,
{\em Indagationes Mathematicae\/} {\bf 34} (1972), 257--262.

\bibitem{W} H.S.~Wilf, The friendship theorem, in: {\em Combinatorial
    Mathematics and its Applications (D.J.A.~Welsh, ed.)\/},
  Academic Press, London, 1971, pp. 307--309.

\end{thebibliography}
\end{document}